\magnification=\magstep1
\def\Z{\bf Z}

\def\s{\Sigma}
\def\N{{\bf N}}

\vskip1cm

\centerline{\bf  Algebraic extensions of an Eilenberg-MacLane
spectrum}

\centerline{\bf in the category of ring spectra }

\medskip

\centerline{\it by Stanislaw Betley}

\vskip1cm

{\bf 1.Introduction.}

\bigskip

The content of the following note lives at the border between
algebra and topology. Historically the origin and  development of
algebraic topology was stimulated by beautiful applications of
algebraic methods for solving topological problems. Later it
turned out that going in opposite direction can be fruitful for
algebra also. Perhaps the first observation of this type can be
derived  from the celebrated  Dold-Kan theorem from the fifties,
which can be viewed as a statement that topological observations
about Eilenberg-MacLane spaces should have meaning in the category
of chain complexes. In proceeding years we observe quick
development of the point of view that topological objects and
methods should give fruitful observations for algebra. We can give
many examples here like various applications in algebra  of
topologically defined homology theories or algebraic $K$-theory,
but of course this is not our aim in this paper.

John Rognes in [R] defined extensions of ring spectra which have
algebraic origin and flavor. So we can talk about Galois
extensions of ring spectra , separable extensions , $thh$-\'etale
extensions and just \'etale ones. When $R$ is a ring we can
associate to it an Eilenberg-Maclane ring spectrum $HR$ so we can
view problems about rings as problems in topology. We would like
to spend some time on studying the following question: do we get
this way any new extensions of an Eilenberg- MacLane spectrum $HR$
for a commutative ring $R$ ? In other words: does every extension
of $HR$ come from an extension of rings of the corresponding type
(Galois, separable, \'etale) ? Speaking again in a different way:
do we get anything new for the theory of rings via embedding them
in the stable homotopy category ?

The following note is mostly devoted to the easy part of the
problem. We are going to show that in the case of Galois
extensions the answer to the question above is negative. So every
Galois extension of $HR$ in the category of spectra comes from
Galois extension of rings. Such a strong statement is not true in
the case of separable extensions. We show that under some
additional assumptions imposed on the extension we get similar
statement. Every ring spectrum comes with an associated graded
ring of homotopy groups. Hence in both cases we first prove the
corresponding statement about graded rings and then about ring
spectra. The graded algebraic case is not necessary for
topological arguments but can serve as a source of some  good
intuitions.

In Section 5 of the paper we approach the case of \'etale type
extensions of spectra. Here we can fully answer our main question
only in the connective case and the answer is the same as for
extensions of  Galois type. On the other hand we know that in
general the situation for etal\'e extensions is different than in
the connective case. In some sense this was the crucial
observation of Mandell,  which was the starting point for the
consideration of this note. As discussed in [MM, example 3.5],
Mandell in private communication showed that the extension
$HF_p\to B$ is \'etale (in certain sense) where $B=F(K(Z/p,n),
HF_p)$ is a mod p cochain $HF_p$-algebra of an Eilenberg-MacLane
space $K(Z/p,n)$ for $n\geq 2$.

In the paper we use freely language of [R] and [EKMM]. So while in
topological world  we work in the category of $S$-algebras and
$S$-modules, where $S$ of course denotes the sphere spectrum. In
algebra all our rings are unital with unital maps.

\bigskip

{\bf Acknowledgment:} This research was partially supported by the
Polish Scientific Grant N N201 387034.

\vskip1cm

{\bf 2. Preliminaries on spectra}.

\bigskip

In this short section we try to put all necessary notation needed
for the rest of the paper. As was said in the introduction our
basic reference is [EKMM]. Let us recall briefly from there what
we mean under the word ''spectrum''. Assume that we have a
structure of a real inner product space on $R^{\infty}$. Then  a
spectrum $E$ is a way of associating a based space to every finite
dimensional vector space $V\subset R^{\infty}$ with a structure
homeomorphisms
$$\sigma_{V,W}:EV\to \Omega^{W-V}EW$$
when $V\subset W$. Here $W-V$ is the orthogonal complement of $V$
in $W$ and $\Omega^WX$ is the space of based maps from $S^W$ to
$X$, where $S^W$ is a one point compactification of $W$. The map
between spectra is just a family of based maps index by $V\subset
R^{\infty}$ commuting with the structure maps. This way we obtain
a category of spectra $\cal S$. The functor from spectra to spaces
given by restriction to the $V$th space has a left adjoint which
is denoted in [EKMM]  by $E^{\infty}_V$, or $E^{\infty}_n$ in the
case when $V=R^n$. When $V=0$ it is a suspension functor
$E^{\infty}$. For $0\leq n$ we define the spectrum $n$-sphere
$S^n$ as $E^{\infty}S^{R^n}$. For $0>n$ we define $S^n$ as
$E^{\infty}_{-n}S^0$. For $m\geq 0$ there are canonical
isomorphisms $E^mS^{R^n}\simeq S^{m+n}$ and $E^{\infty}_mS^n\simeq
S^{n-m}$. We define homotopy groups of a spectrum $E$ to be
$$\pi_n(E)=h{\cal S}(S^n,E)$$
where $h$ stands here for the homotopy category of spectra,
obtained via correct choice of a closed model category structure
on ${\cal S}$.

The letter $S$ denotes the $0$-sphere spectrum, as defined in
[EKMM, Section 1], called usually just as sphere spectrum.  It
comes with a map $S\wedge S\to S$ giving it a structure of a ring
spectrum. Throughout the paper we  are working in the category of
$S$-algebras, as defined in [EKMM, Section 2 and 3]. It means that
we are considering $S$-modules, denoted by capital letters $A$,
$B$, etc.,  equipped additionally with the multiplication and unit
maps

$$\mu: A\wedge_S A\to A$$

$$1_A:S\to A$$
satisfying standard associativity and unity conditions. When $A$
is an   $S$-algebra  we can define the symmetric monoidal category
${\cal M}_A$ of right $A$-modules (left $A$-modules $_A{\cal M})$.
As objects in it we have $S$-algebras equipped additionally with
the right action of $A$ on them:
$$B\wedge_S A\to B$$ satisfying standard associativity and unity
conditions, which we know from the algebraic category of modules.
We say that $B$ is an $A$-algebra if it is a monoid in ${\cal
M}_A$ (compare [EKMM, Section 7]).
 In case $A$ is
not commutative we can also talk about categories of right and
left $A$-modules but  our algebras will be always unital by what
we mean that they come equipped with a unit map

$$1_B^A:A\to B$$
which is compatible with $1_A$ and $1_B$. We consistently remove
$S$ from our notation, hence for example $1_A$ is the same as
$1_A^S$, $A\wedge A$ denotes $A\wedge_S A$, etc., etc.

\bigskip

 {\bf 3. Galois  extensions}.

\bigskip
Let $A$ and $B$ are commutative rings and $G$ a finite group.
Following [G] we define:
\bigskip

{\bf Definition 3.1:}  We say that the extension of commutative
rings $A\hookrightarrow B$ is $G$-Galois if $G$ is a subgroup of
$Aut(B/A)$, $B^G=A$ and the map $h:B\otimes_A B\longrightarrow
Map(G,B)$ is a $B$-algebra isomorphism, where $h(x\otimes
y)(g)=x\cdot g(y)$.

\bigskip

In the case of graded rings we assume that the action of $G$
preserves grading. By a grading we always mean here $Z$-grading.
The $B$-algebra of functions $Map(G,B)$ will be also viewed very
often as $\prod_{g\in G} B$ so we can project from it on the
coordinate corresponding to the given $g\in G$. Observe that the
map $h$ preserve natural gradings of the source and the target.
\bigskip

{\bf Theorem 3.2:} Let $A\hookrightarrow B$ be a Galois extension
of graded rings. If $A$ is nontrivial only in grade $0$ then the
same is true for $B$.
\medskip

Proof. We show  first (after [G, Theorem 1.6]) that $B$ is a
finitely generated projective $A$-module. The proof given there
works as well in the graded case. We present it here because the
careful  looking at the formulas from topological point of view
gives us the desired result for spectra. Let $\Sigma x_i\otimes
y_i$ be the preimage of $(1,0,...,0)\in \prod_{g\in G} B$, where
$1$ is at the coordinate corresponding to the unit $e$ of $G$.
Define the $A$-linear trace $tr:B\to A$ by $tr(y)=\Sigma_{g\in
G}g(y)$. Let $\varphi_i:B\to A$ be defined by
$\varphi_i(z)=tr(zy_i)$. Then the direct calculation gives us the
formula for any $z\in B$:
\medskip

(3.2.1)
$$z=\mathop\Sigma_i\varphi_i(z)\cdot x_i$$ This immediately implies that
$B$ is a finitely generated projective $A$-module because formula
(3.2.1) shows that the pairs $(x_i, \varphi_i)$ form a dual basis
for $B$ over $A$. But for us it is  more important to observe that
formula (3.2.1) shows that $B$ can have non trivial elements only
in finitely many gradations (is finitely graded) because we have
only finitely many $x_i$s and $A$ is fully in the $0$-grade. This
observation immediately implies our statement. If $k$ is the
highest (lowest) nontrivial gradation  of $B$ then $B\otimes_A B$
has highest (lowest) nontrivial gradation in dimension $2k$. But
grading of $\prod_{g\in G} B$ is the same as the grading of $B$.
Hence $k$ has to be $0$.

\bigskip
Now we move towards topology. Let $A\to B$ be a map of commutative
$S$-algebras and $G$ is a finite group acting continuously from
the left on $B$ via the $A$-algebra maps.  Let us recall (compare
 [R, Definition 4.1.3]) the definition of the  Galois extensions in the
category of $S$-algebras.

\medskip

{\bf Definition 3.3:} With the assumptions as above we say that
$A\to B$ is a $G$-Galois extension of $S$-algebras if two
canonical maps of $S$-modules $i:A\to B^{hG}$ and $h:B\wedge_A
B\to F(G_+,B)$ are weak equivalences.

\medskip
Perhaps we should also recall here after [R] the definitions of
the maps $i$ and $h$. The map $i:A\to B^{hG}=F(EG_+,B)^G$ is the
right adjoint to the composite $G$-equivariant map $A\wedge EG_+
\to A\to B$, collapsing the contractible free G-space to a point.
The map $h$ is right adjoint to the composite map $B\wedge_A
B\wedge G_+ \to B\wedge_A B \to B$ where the first map come from
the action of $G$ on the middle $B$ from the left and the second
is just the multiplication map. Observe that in our case ($G$
finite) we can equally well write $F(G_+,B)$ as $\prod_{g\in G}
B$. Note also that we can view $h$ as
$$B\wedge_A B\ \ \buildrel{id\wedge \prod {g}}\over \longrightarrow
\ \ B\wedge_A \prod_{g\in G} B \ \buildrel{\prod\mu \circ
(id\wedge pr_g)} \over \longrightarrow\  \mathop\prod_{g\in G} B
$$

\noindent where we denote by $g$ the map $B\to B$ coming from the
action of $g\in G$ on $B$ and $pr_g$ denotes the projection on the
$g$-factor.

\bigskip
{\bf Theorem 3.4:} Let $A=HR\to B$ be a $G$-Galois extension of
commutative ring spectra. Then $B$ is equivalent to $H(\pi_0B)$
and $R\to \pi_0B$ is a $G$-Galois extension of commutative rings.

\medskip
Proof. The proof is a combination of results from [R] and [EKMM].
By [R, Proposition 6.2.1] we know that $B$ is a dualizable
$A$-module. Then by [R, Proposition 3.3.3] and [EKMM Chapter III,
Theorem 7.9] we know that $B$ is a retract of a finite cell
$A$-module. This implies that $B$ has only finitely many
nontrivial homotopy groups each of which is a finitely generated
$R$-module.

We will prove first that $B$ is an Eilenberg-MacLane spectrum. We
will follow the lines of the algebraic  graded case,  the argument
is only a little more delicate. On the other hand this is the
crucial step because the rest of our theorem is then proved in [R,
Theorem 4.2.1]. For the readers convenience we will sketch Rognes'
argument later.

Let $k$ be the lowest integer such that $\pi_k(B)\neq 0$. Assume
that $k<0$. Then by the spectral sequence for the homotopy groups
of a smash product, which is described below, we know that in
$\pi_{2k}({B\wedge_A B})$ we have classes coming from $
\pi_k(B)\otimes_R \pi_k(B)$. If this latter group is nontrivial we
get a contradiction as in the graded algebraic case. But the group
in question indeed is nontrivial  by a simple algebraic lemma,
probably well known to everybody:

\bigskip

{\bf Lemma 3.4.1:} Assume that $T$ is a commutative ring and $M$
is a finitely generated module over $T$. Then $M\otimes_T M$ is
nontrivial.

\medskip

Proof. Assume that $M$ has only one generator. Then $M$ is
isomorphic to $R/I$ for a certain ideal $I$. Let $J$ be a maximal
ideal containing $I$ then $R/I$ maps epimorphically onto $R/J$. We
know that $R/J\otimes_T R/J$ is nontrivial by maximality of $J$
(is isomorphic to $R/J$) so by right exactness of the tensor
product we know that $R/I\otimes_T R/J$ is nontrivial. Hence,
again by the right-exactness of the tensor product we get that
$R/I\otimes_T R/I$ is nontrivial.

We can proceed further by induction with respect to the number of
generators in $M$. If $M$ has $n$ generators then it fits into an
exact sequence of $T$-modules

$$0\to L\to M\to N\to 0$$

\noindent in which $L$ has one and $N$ has $n-1$ generators. By
induction $N\otimes_T N$ is nontrivial and $M\otimes_T M$ maps
epimorphically onto $M\otimes_T N$ which maps onto $N\otimes_T N$
by the left-exactness of the tensor product. So the proof of our
lemma is finished.

\bigskip

Now we come back to the proof of 3.4. If $k\geq 0$, and  hence $B$
is connective, we know by [EKMM IV, Proposition 1.4] that the dual
$A$-spectrum of $B$ is coconnective (has nontrivial homotopy
groups only in nonpositive dimensions). On the other hand by [R,
Proposition 6.4.7] $B$ is self dual. So homotopy groups of $B$
have to be concentrated in dimension $0$, as we wanted to show.

Now we can finish the proof of 3.4. Since we know now that $B$ is
an Eilenberg-MacLane spectrum we can recall [R,Proposition 4.2.1].
Let us write $T$ for $\pi_0(B)$ for shortness. By [EKMM, IV.4.3]
we have the homotopy fixed point spectral sequence

$$E^2_{s,t}=H^{-s}(G,\pi_tHT)\Longrightarrow \pi_{s+t}(HT^{hG})$$
which in our case gives us $T^G\simeq \pi_0(HT^{hG})\simeq
\pi_0(HR)=R$

Similarly we have useful  spectral sequence for the homotopy
groups of a smash product, which was used before and will be
crucial in the next section. It is of the form

$$E^2_{s,t}=Tor^R_{s,t}(T,T)\Longrightarrow \pi_{s+t}(HT\wedge_{HR} HT)$$

It gives $T\otimes_RT\simeq \pi_0(HT\wedge_{HR} HT)\simeq
\pi_0(\prod_{g\in G} HT)=\prod_{g\in G} T$. This implies that
$R\to T$ is $G$-Galois in the algebraic sense.

\bigskip

{\bf Remark 3.5:} The proof that $B_*$  is finitely generated over
$R$ is not direct. We would like to present below a sketch of a
direct argument which mimics the proof of the algebraic equivalent
statement used in the proof of 3.2.

Our extension is $G$-Galois so the map $h$ defined before is a
weak equivalence. This means that the unit map $1_B:S\to B$ can be
factored as
\bigskip

3.5.1

$$S\buildrel{\varphi} \over\longrightarrow B\wedge_A B \buildrel{h}
\over\longrightarrow \prod_{g\in G} B\buildrel{pr_e}
\over\longrightarrow B$$

\noindent where $\varphi$ is the preimage of $(1,0,...,0)\in
\pi_0(\prod_{g\in G} B)$ and, as was defined before,  $pr_e$ is
the projection map on the coordinate corresponding to the trivial
element $e\in G$. Equivalently, by the choice of $\varphi$,  we
could say that $1_B$ can be factored as

\bigskip

3.5.2

$$S\buildrel{\varphi} \over\longrightarrow B\wedge_A B \buildrel{h}
\over\longrightarrow \prod_{g\in G} B\buildrel{\oplus 1_B}
\over\longrightarrow B$$

Let $f:S^n\to B$ be a map representing an element in $\pi_n(B)$.
Then
$$S\wedge S^n\buildrel{1_B\wedge f}\over \longrightarrow B\wedge_A
B\buildrel{\mu}\over \longrightarrow B$$ represents the same
element in $\pi_n(B)$ as $f$. But instead of $1_B$ we can use the
composition of maps from 3.5.1 or 3.5.2. The composition $
B\buildrel{\Delta }\over \longrightarrow \prod_{g\in G}
B\buildrel{\oplus g}\over \longrightarrow B$  will be denoted by
$\phi$ in the future. Observe that the following two maps:

\bigskip

3.5.3

$$B\wedge_A B\wedge_A B\buildrel{h\wedge id}\over \longrightarrow
(\prod_{g\in G} B)\wedge_A B\buildrel{(\oplus id)\wedge id}\over
\longrightarrow B\wedge_A B\buildrel{\mu}\over\longrightarrow B$$

\noindent and

\bigskip

3.5.4

$$B\wedge_A B\wedge_A B\buildrel{id\wedge \mu}\over \longrightarrow
 B\wedge_A B\buildrel{id \wedge \phi}\over \longrightarrow
B\wedge_A B\buildrel{\mu}\over\longrightarrow B$$

\noindent are homotopic after precomposing with $\varphi\wedge f$.
This follows immediately from the definition of $\varphi$. The map
3.5.3 precomposed with $\varphi\wedge f$ is homotopic to $f$. On
the other hand the map $\phi $, as being $G$-invariant, factors
through the spectrum $B^{hG}$ which is equivalent to $A$. Hence
the homotopy properties of 3.5.4 precomposed with $\varphi \wedge
f$ depend only on the homotopy class of $\varphi$ and homotopy
groups of A. This implies immediately that $B$ can have only
finitely many nontrivial homotopy groups.

\bigskip

{\bf 4. Separable extensions}.

\bigskip

For separable extensions of ring spectra we would like to prove
the same statement as was proved for  Galois extensions in the
previous section. From the ideological point of view this is the
expected statement because as in algebra one expects that any
commutative separable extension embeds into a $G$-Galois one, for
a certain $G$. We  are able to get the  expected result under
additional hypothesis on the  extension. We expect that this
result is well known to experts but we could not find any place
with a proof written down. But before going into stable homotopy
category let us state and prove the graded algebraic counterpart
of this statement. Later we will generalize the proof to the case
of spectra. Let $A$ and $B$ be two $\Z$-graded unital rings.

\medskip

{\bf Definition 4.1:} We say that $A\to B$ is separable if the
$A$-algebra multiplication map $\mu:B\otimes_A B^{op}\to B$,
considered as a map in the  category of $B$-bimodules, admits a
section $\sigma:B\to B\otimes_A B^{op}$.
\bigskip

{\bf Theorem 4.2:} Assume that $A$ is concentrated in gradation
$0$ only. Let $A\to B$ be a separable extension of graded rings as
defined above and $B$ has no zero divisors in the subring $B_0$.
Then $B$ is concentrated in gradation $0$.

\medskip
Proof. The crucial but obvious observation in the case $A=A_0$ is
 that if $x_1\otimes x_2=x_3\otimes x_4 \neq 0$ in $B\otimes_A B$ and all
$x_i's$ are of homogeneous degree then $deg(x_1)=deg(x_3)$ and
$deg(x_2)=deg(x_4)$. This is the case because $B\otimes _A B$ has
double grading and multiplication by elements of $A$ preserves it.
Separability means that there exists an element
$$\mathop\Sigma_{i=1}^kb_i\otimes c_i\ \in \ (B\otimes B)_0 \subset B\otimes_A B$$
satisfying for any $b\in B$

$$b(\mathop\Sigma_{i=1}^kb_i\otimes c_i)=(\mathop\Sigma_{i=1}^kb_i\otimes
c_i)b$$ The element described above is equal to the image of $1$
under the map $\sigma: B\to B\otimes_A B$.  We can assume that the
elements $b_i$ and $c_i$ are homogeneous and $deg(c_i)=-deg(b_i)$.
Let  $\{ b_{i_j}\}_{j\in J}$ be the set of these $b_is$ which have
the highest grade. Then $b(\s_{j\in J}b_{i_j}\otimes c_{i_j})=0$
for any $b$ of grade bigger than $0$ by degree reasons. The same
one can say about any $b$ of negative degree considering
$(\s_{j\in J}b_{i_j}\otimes c_{i_j})b=0$ . Hence either $b$ has to
be zero  or by our assumption on $B_0$, $\Sigma_{j\in
J}b_{i_j}\cdot c_{i_j}=0$ and we can send $1$ to
$$\mathop\Sigma_{i=1}^kb_i\otimes c_i-\mathop\Sigma_{j\in J}b_{i_j}\otimes
c_{i_j}$$

But this latter element has the lower highest degree among $b_is$
so step by step we can lower this highest degree in our sum to
$0$. Obviously then all $c_is$ should have also degree $0$. This
means, we can assume that $\sigma(1)\in B_0\otimes B_0$. But then
in order to have satisfied

$$b(\mathop\Sigma_{i=1}^kb_i\otimes c_i)=(\mathop\Sigma_{i=1}^kb_i\otimes
c_i)b$$

an element $b$ cannot have degree different from $0$.
 Thus $B$ should have only
$0$ grade and $A_0\to B_0$ should be a separable extension of
ungraded rings.

\medskip

{\bf Remark 4.3:} It is easy to observe that if our rings are only
$\N$-graded (connective) then  theorem 4.2 is true without any
assumption on $B_0$.

\bigskip
Now we come to the  definition of separable extension of ring
spectra, as it is given in [R, Definition 9.1.1].

\medskip

{\bf Definition 4.4:} We say that $A\to B$ is separable if the
$A$-algebra multiplication map $\mu:B\wedge_A B^{op}\to B$,
considered  as a map in the stable homotopy category of
$B$-bimodules relative to $A$, admits a section $\sigma:B\to
B\wedge_A B^{op}$.

\bigskip

Observe that in our case if $B$ is an extension of $HR$ then being
a module over an Eilenberg-MacLane spectrum it is equivalent  in
the stable homotopy category to the wedge of Eilenberg-Maclane
spectra $H(B_i,i)$. Of course here $B_i=\pi_i(B)$, every $B_i$
carries a structure of an $R$-module and hence $H(B_i,i)$ carries
a structure of $HR$-module as well. Let $v:\bigvee H(B_i,i) \to B$
gives us an equivalence guaranteed above.

\medskip

{\bf Theorem 4.5:} Let $HR\to B$ be a separable extension as
defined above. Assume, similarly as previously, that $\pi_*(B)$
has no $0$-divisors in $\pi_0(B)$.  Assume moreover that the map
$v$ described above is an $HR$-module map.  Then $B$ is equivalent
to $H\pi_0B$ and $R\to \pi_0B$ is a separable extension of rings.

\medskip

Proof. We would like to follow the lines of the proof of 4.2
taking as an extension of $R$ the ring $B_*=\pi_*(B)$. The map
$\sigma$ gives us the splitting of $\pi_*(B)$ from
$\pi_*(B\wedge_A B^{op})$.  But we are not able to use 4.2
directly because the statement $\pi_*(B\wedge_A
B^{op})=\pi_*(B)\otimes_{\pi_*(A)} \pi_*(B)$ is false in general.
Instead, as it was mentioned in the previous section, we have only
a spectral sequence converging to $\pi_*(B\wedge_A B^{op})$
([EKMM, chapter IV] with the second table given by the formula

$$E_{p,q}^2=Tor_p^{\pi_*(A)}(B_*,B_*)_q$$

In order to apply similar argument as previously we are forced to
study the bimodule (kind of) structure of this second table over
$B_*$. The main point is that $\sigma$ is a bimodule map so for
every element $b\in \pi_*(B)$ we have as previously
$$\sigma_*(b)=b\sigma_*(1)=\sigma_*(1)b$$

\noindent Hence the multiplication by $b$  on $\sigma_*(1)\in
\pi_*(B\wedge_A B)$ should have the same effect when we use right
and left module structures. Our ground ring spectrum is $HR$ hence
every group $B_i$ is a module over $R$ and graded $R$-resolution
of $B_*$ is just a graded sum of ordinary $R$-resolutions of
$B_i$s. This leads to the splitting formula for $Tor$-groups:

$$Tor_p^R(B_*,B_*)=\bigoplus_{i,j}Tor_p^R(B_i,B_j)$$
and
$$E^2_{p,q}=\bigoplus_{i+j=q}Tor_p^R(B_i,B_j)$$
Moreover, because of our hypothesis on the map $v$, $B$  is
homotopically a wedge of Eilenberg-MacLane spectra itself. Hence
we can write a resolution of $B_*$ coming from the wedge of
resolutions of $H(B_i,i)$s. This leads to the resolution of $B_*$
which is a sum of resolutions of $B_i$s.  Then the spectral
sequence above can be viewed as a sum of spectral sequences coming
from $\pi_*(S^i\wedge HB_i\wedge B)$. This implies also that the
limit group splits into a corresponding sum.

\noindent Let $b$ be represented by a map $S^n\to B$, which will
be called also $b$. Then the multiplication by $b$ on $B_*$ from
the left  can be described by a map
$$S^n\wedge B\to B\wedge B\to B\wedge_{HR} B\to B$$
where the first map is $b\wedge id$, the next comes from the fact
that $S$-module and $HR$-module structures on $B$ are compatible,
the last map is $\mu$ and $\wedge$ without subscript denotes
smashing over the sphere spectrum. Of course we have similar
description for the right multiplication. We can smash the
sequence above with $B$ over $HR$  and get a natural map
$S^n\wedge B\wedge_{HR} B \to \ B\wedge_{HR} B$. Similarly we have
a natural map $B\wedge_{HR} B\wedge S^n \to \ B\wedge_{HR} B$ and
these two maps give us two maps of spectral sequences (for  the
functoriality of the spectral sequence construction see [EKMM,
Section IV.5])
$$_lE^2_{p,q}=Tor^{R}_p(\pi_*(S^n\wedge B),\pi_*(B))_q \to Tor^{R}_p(\pi_*(B),\pi_*(B))_q=E^2_{p,q}$$
and
$$_rE^2_{p,q}=Tor^{R}_p(\pi_*(B),\pi_*(B\wedge S^n))_q \to Tor^{R}_p(\pi_*(B),\pi_*(B))_q=E^2_{p,q}$$
\noindent where the letters $l$ and $r$ refer to the left and
right multiplication by $b$. Observe that $\pi_*(S^n\wedge B\wedge
B)=\pi_{*-n}(B\wedge B)$ as one is an n-fold suspension of
another. Spectral sequences $_lE^2_{p,q}$ and $_rE^2_{p,q}$ are
the same as $E^2_{p,q}$ with a shift of total grading by $n$. The
described above map $_lE^2_{p,q}\to E^2_{p,q}$ is induced by
multiplication with $b$ and on the level of the algebraic second
table of  spectral sequences it is induced by the map $\cdot
b:B_{*-n}\to B_*$ on the first variable in the groups $Tor$. One
has the same description for the map $_rE^2_{p,q}\to E^2_{p,q}$
but the multiplication goes along the second variable in the
$Tor$-groups. So starting from $Tor_p(B_i,B_j)$ we go by
multiplication by $b$ one time to $Tor_p(B_{i+n},B_j)$ and the
second to $Tor_p(B_i,B_{j+n})$. When we know all of this we can
apply the same procedure as in 4.2 to the image of $1\in \pi_0(B)$
in $\pi_*(B\wedge_{HR} B)$.

 Multiplication from the left by
an element $b\in B_n$ takes $B_i$ to $B_{i+n}$ and this is a map
of right $R$-modules hence extends to the map of resolutions. It
means that multiplication with $b$ from the left induces a map of
spectral sequences

$$^iE^2_{p,q}\to ^{i+n}E^2_{p,q}$$

On the other hand multiplication from the right by $b$ is a map of
coefficients which is a map of left $R$ modules. It means it
induces self maps of $^iE^2_{p,q}$'s. Now we can easily argue that
if $b$ has grade different from $0$ then multiplication by it on
$\sigma_*(1)\in \pi_*(B\wedge_{HR} B)$ should be trivial. As was
stated before $\sigma_*(1)$ can be decomposed into a finite sum of
nontrivial elements $\Sigma_{j=1}^k x_{i_j}$ , where each piece
$x_{i_j}$ comes from $^{i_j}E^2_{p,q}$ for different $i_j$'s. Let
$k$ be the highest among $i_j$'s. Then $b\cdot \sigma_*(1)$ has a
summand coming from $^{k+n}E^2_{p,q}$. On the other hand
$\sigma_*(1)\cdot b$ does not have such a summand. It means that
$b\cdot x_k$ should be $0$ for any $b$ of degree higher than $0$.
Observe that by definition $x_k$ is of homotopical degree $0$ so
by our hypothesis on $B_*$ if $b$ is non zero then $x_k$ must be
$0$.This way we can lower the maximal index $k$ in the
decomposition of $\sigma_*(1)$. But this argument works always
when $B_*$ has elements of degree higher that $0$ so we know that
$B_*$ should be concentrated in non-positive degrees. But of
course for $b$ of negative degree we argue similarly, starting
from $x_m$ with the lowest possible index $m$. This implies that
$B$ is an Eilenberg-MacLane spectrum and the extension $HR\to B$
comes from a separable extension of $R\to R'$ where obviously
$B\simeq HR'$.

\bigskip

{\bf Example 4.6}(after Birgit Richter): Let $B = F_{HR}(\Sigma HR
\vee \Sigma^{-1} HR, \Sigma HR \vee \Sigma^{-1} HR)$, where
$F(.,.)$ denotes the mapping spectrum in appropriate category.
This is a Brauer-trivial Azumaya algebra over $HR$, so in
particular it is separable over $HR$. Observe that
$B_*(=\pi_*(B))$ is isomorphic to $R$ for $*=2, -2$ and
$B_0=R\oplus R$. One checks directly that $B_*$ is isomorphic as
an $R$-algebra to the algebra of $2\times 2$-matrices over $R$
with appropriate grading. So $B_*$ is separable over $R$. Hence it
is difficult to imagine how one could make  our assumption on
$B_0$ in 4.2 or 4.5 weaker.

\bigskip

{\bf 5. \'Etale extensions}.

\bigskip

Let us start this section from the definition of the topological
Hochschild homology(see [R, section 9.2] or [EKMM, Chapter IX]):

\bigskip
{\bf Definition 5.1:} Let $B$ be an algebra over a commutative
$S$-algebra $A$. Then we define
$$THH^A(B)=Tor^{B\wedge_A B^{op}}(B,B)$$

\medskip

{\bf Definition 5.2:} We say that $A\to B$ is formally
symmetrically \'etale if the canonical map $\zeta :B\to THH^A(B)$
is a weak equivalence.

\bigskip

For an extension $A\to B$  of ordinary rings we should use the
same definition of topological Hochschild homology as above for an
extension $HA\to HB$. This is because $THH$-theory does not have
algebraic definition. But of course this leads  to problems with
the graded algebraic case because we do not know what is $HB$ for
a graded ring $B$. Define $THH^A_i(B)=\pi_i(THH^A(B))$, which has
precise meaning in the ungraded case . Then $THH_1^A(B)=HH_1(B)$
where $HH$ denotes the ordinary Hochschild homology. Hochschild
homology theory has perfect meaning also in the graded case and so
we  can say something about the graded rings using Hochschild
homology groups. Assume that $B$ is a graded ring which is an
algebra over a commutative ring $A$. We treat $A$ as a graded
object concentrated in gradation $0$. Moreover assume
 that $B$ is commutative and trivial in negative gradations. Then

\bigskip

{\bf Lemma 5.3:} If $B$ is nontrivial in positive gradations and
$B_0=A$ then $HH_1(B)$ is nontrivial. Hence we can think about $B$
as being  not formally symmetrically \'etale over $A$.

\medskip

Proof: The first  Hochschild homology group of $B$ is the same as
 the group of K$\ddot a$hler differentials of $B$ over $A$. If
$B_k$ is the lowest positive nontrivial gradation of $B$ then
$(\Omega^1_{B/A})\supset B_k \neq 0$. This follows directly from
the definition of K$\ddot a$hler differentials. If $b\in B_k$ is
of the form $b_1\cdot b_2$ then either $b_1\in A$ or $b_2\in A$ by
our choice of $k$. Hence there are no relations except linearity
over $A$ between generators $db$ of $(\Omega^1_{B/A})$ for the
elements $b\in B_k$.
\bigskip

The lemma above shows the way in which we can approach the similar
problem for ring spectra. The key ingredient is hidden in the
spectrum  $\Omega_{B/A}$ of differential forms of $B$ over $A$
defined as a cofibrant replacement of the homotopy fiber of the
multiplication map $\mu: B\wedge_A B^{op}\to B$. Assume that $B$
is a connective $HA$-algebra  where $HA$ is an Eilenberg-MacLane
spectrum of a commutative ring $A$. Assume that $\pi_0(B)=A$. Then

\bigskip

{\bf Theorem 5.4:}  If $B$ has higher nontrivial homotopy groups
then $B$ is not a formally  symmetrically \'etale extension of
$HA$.

\medskip
Proof. Let $k$ be the smallest natural number bigger than $0$ for
which $\pi_k(B)\neq 0$. Then from the spectral sequence for
$\pi_*(B\wedge_{HA} B)$ we know that $\pi_k(B\wedge_{HA}
B)=\pi_k(B)\oplus \pi_k(B)$. Observe that there is a map

$$i:B\to B\wedge_{HA} B$$ which composed with multiplication
$$\mu :B\wedge_{HA} B\to B$$
is trivial. It is the topological counterpart of the algebraic map
$X\to X\otimes X$ which takes $x\in X$ to $x\otimes 1-1\otimes x$.
The map $i$ is defined as the difference  of maps $id_B\wedge 1_B$
and $1_B\wedge id_B$. Of course $\mu\circ i=0$ and hence $i$
factors as $j\circ \beta$ through the homotopy fiber
$j:\Omega_{B/A}\to B\wedge_{HA} B$ of $\mu$. Let $\alpha:S^k\to B$
represents the nontrivial element in $\pi_k(B)$. Then $i\circ
\alpha$ is nontrivial on homotopy groups by the formula for
$\pi_k(B\wedge_{HA} B)$ and hence $\beta \circ \alpha$ is also
nontrivial. It means that $\Omega_{B/A}$ is not contractible and
has nontrivial $k$-th homotopy group. We know that $B$,
$B\wedge_{HA} B^{op}$ and  $\Omega_{B/A}$ are connective and it
follows that the latter spectrum is $(k-1)$-connected from the
definition of $\Omega_{B/A}$. Moreover we know that
$\pi_0(B)=A=\pi_0(B\wedge_{HA} B^{op})$ and $\pi_i(B\wedge_{HA}
B^{op})=0$ for $i=1,...,k-1$. This last calculation follows
directly from the spectral sequence for calculating homotopy
groups of $B\wedge_{HA} B^{op}$.

Again , by the definition of $\Omega_{B/A}$ we have a cofiber
sequence

$$B\wedge_{B\wedge_{HA} B^{op}} \Omega_{B/A} \to B \buildrel{\zeta} \over\longrightarrow THH^{HA}(B) $$

\noindent Our proof will be finished if we show that
$B\wedge_{B\wedge_{HA} B^{op}} \Omega_{B/A} $ is not weakly
equivalent to a point. But again we can use the spectral sequence
for calculating  homotopy groups of this spectrum with the second
table given by the formula

$$E^2_{s,t}=Tor_{s,t}^{\pi_*(B\wedge_{HA} B^{op})}(B,\Omega_{B/A})$$

\noindent Taking into account the connectivity of $\Omega_{B/A}$
and assumptions on $B$ we immediately read that
$\pi_k(B\wedge_{B\wedge_{HA} B^{op}}\Omega_{B/A})$ equals to the
$k$th grade of $\pi_*(B)\otimes_{\pi_*(B\wedge_{HA} B^{op})}
\pi_*(\Omega_{B/A})$. But this latter group is easily calculated
by dimension reasons as $A\otimes_A \pi_k(\Omega_{B/A}) =
\pi_k(\Omega_{B/A})\neq 0$.

\bigskip

{\bf Remark 5.5:} Observe that in both Lemma 5.3 and Theorem 5.4
above the assumption on the $0$-grade is irrelevant. In 5.3 by the
described above arguments one gets the  result for differential
forms of $B$ over $B_0$. But then either $0\neq
(\Omega^1_{B_0/A})\subset (\Omega^1_{B/A})$  or
$(\Omega^1_{B_0/A})=0$ and then by the same argument as previously
we get $(\Omega^1_{B/A})\supset B_k \neq 0$. The argument for 5.4
is left to the interested  reader.

\bigskip

{\bf Remark 5.6:} In the commutative case Rognes defined the
notion of formally \'etale extension $A\to B$ using the notion of
topological Andr\'e-Quillen homology. We do not recall it here
because the property of being formally symmetrically \'etale is
equivalent to formally \'etale for connective algebras. So far we
are not able to analyze the non-connective cases, where we know
that these two notions are different by Mandell's example.

\bigskip

{\bf Cojecture 5.7:} Theorem 5.4 is true for $HA$-algebras which
are bounded below.

\medskip

As an evidence we we show below that we can extend the proof of
the lemma 5.3 to cover algebras which are bounded below. We have
the following lemma:

\bigskip

{\bf Lemma 5.8:} Assume that $B$ is a graded commutative
$A$-algebra, $B_0=A$ and $B_i=0$ for $i<k$ where $k$ is some
negative number. Then $HH_1(B)$ is nontrivial.

\medskip

Proof: As previously we will use the fact that $HH_1(B)$ is equal
to the $B$-module of K\"ahler differentials and the latter module
is the same as $I/I^2$ where $I$ is a kernel of the multiplication
map $B\otimes_A B\to B$. Assume that $x\in B_k$ and $dx=x\otimes
1-1\otimes x$ is in $I^2$. Then by the argument from the beginning
of the proof of 4.2 and the fact that $k$ is the minimal grade
with nontrivial $B_k$ we get immediately that $x$ can be expressed
as an $A$-combination of elements of the form $x_1\cdot x_2$ where
both $x_1$ and $x_2$ have negative grading. So if for every $x\in
B_k$ the differential $dx$ is trivial in $HH_1(B)$ then all
elements of $B_k$ are sums of multiples of elements of higher but
negative degree. We can extend this reasoning easily to other
negative degrees of $B$ by induction and get that if the
differential $dx=0$ for $x\in B_s$, $s$ negative,  then $x$ can be
expressed as a sum of finite multiples of elements of negative but
higher than $s$ degree. So either we have nontrivial elements of
negative degree in $HH_1(B)$ or negative part of $B$ is generated
over $A$ by $B_{t}$ where $B_t$ is the highest lower than $0$
nontrivial grade of $B$. But if this is the case then the only
relations among $dx$ for $x\in B_t$ are relations of $A$
linearity. Hence we have nontrivial elements in degree $t$ of
$HH_1(B)$.

\bigskip

{\bf Remark 5.9:} Assume that $B$ is a bounded below $HA$-algebra
with $\pi_0(B)=A$ and for any $i$ , $\pi_i(B)$ is a projective
$A$-module. Then one can easily apply the way of reasoning from
the proof of 5.8 to showing that theorem 5.4 can be extended to
cover such a case.

\bigskip
\vfill\eject

{\bf References:}

\medskip

[EKMM] A.D. Elmendorf, I. Kriz, M.A. Mandell, J.P. May. {\it
Rings, Modules and Algebras in Stable Homotopy Theory}.
Mathematical Surveys and Monographs 47, AMS 1997.

[G] C. Greither. {\it Cyclic Galois extensions of commutative
rings}. LNM 1534, Berlin, Springer-Verlag, 1992.

[MM] R. McCarthy, V. Minasian. {\it HKR theorem for smooth
$S$-algebras}. J. of Pure and Applied Algebra 185 (2003) 239-258.

[R] J. Rognes. {\it Galois extensions of Structured Ring Spectra}.
Mem. AMS 192 (2008).

\bigskip

\bf Instytut Matematyki, University of Warsaw

ul. Banacha 2, 02-097 Warsaw, Poland

e-mail: betley@mimuw.edu.pl

\end